\documentclass[12pt]{article}
\usepackage{amssymb, amsmath, amsfonts,dsfont, mathrsfs,euscript,eufrak,colonequals,indentfirst}
\usepackage{textcomp}
\usepackage{latexsym}

\newtheorem{Theorem}{Theorem}
\newtheorem{Lemma}{Lemma}
\newtheorem{Proposition}{Proposition}
\newtheorem{Remark}{Remark}

\newcommand{\N}{\mathbb N}
\newcommand{\R}{\mathbb R}

\newcommand{\Expec}{ \mathds{E}\,}

\newcommand{\defeq}{\colonequals}

\newcommand{\dd}{\mathrm{d}}
\newcommand{\id}{\mathds{1}}

\newcommand{\AppClass}{\mathcal{A}}
\newcommand{\CovFunc}{\mathcal{K}}
\newcommand{\ContSet}{\mathcal{C}}
\newcommand{\Probab}{\mathds{P}}
\newcommand{\Var}{\mathrm{Var}}

\newcommand{\e}{\varepsilon}
\newcommand{\s}{\sigma}
\newcommand{\cen}{c}
\newcommand{\var}{v}

\begin{document}
\author{A. A. Khartov\footnote{Smolensk State University, 4 Przhevalsky st., 214000 Smolensk, Russia, e-mail: \texttt{alexeykhartov@gmail.com}}\qquad I. A. Limar\footnote{Saint-Petersburg National Research University of Information Technologies, Mechanics and Optics (ITMO University), 49 Kronverksky Pr., 197101 Saint-Petersburg, Russia, e-mail: \texttt{ivan.limar95@gmail.com} }}

\title{Asymptotic analysis in multivariate average case  approximation with Gaussian kernels}
\maketitle

\begin{abstract}
	We consider tensor product random fields $Y_d$, $d\in\N$, whose covariance functions are Gaussian kernels with a given sequence of length scale parameters. The average case approximation complexity $n^{Y_d}(\e)$ is defined as the  minimal number of evaluations of arbitrary linear functionals needed to approximate $Y_d$, with relative $2$-average error not exceeding a given threshold $\e\in(0,1)$. We investigate the growth of $n^{Y_d}(\e)$ for arbitrary fixed $\e\in(0,1)$ and $d\to\infty$. 
	Namely, we find criteria of boundedness for $n^{Y_d}(\e)$ on $d$ and of tending $n^{Y_d}(\e)\to\infty$, $d\to\infty$, for any fixed $\e\in(0,1)$. In the latter case we obtain necessary and sufficient conditions for the following logarithmic asymptotics
	\begin{eqnarray*}
		\ln n^{Y_d}(\e)= a_d+q(\e)b_d+o(b_d),\quad d\to\infty,
	\end{eqnarray*}
	with any $\e\in(0,1)$. Here $q\colon (0,1)\to\R$ is a non-decreasing function, $(a_d)_{d\in\N}$ is a sequence and $(b_d)_{d\in\N}$ is a positive sequence such that $b_d\to\infty$, $d\to\infty$.
	We show that only special quantiles of self-decomposable  distribution functions appear as functions $q$ in a given asymptotics. These general results apply to $n^{Y_d}(\e)$ under particular assumptions on the length scale parameters.
\end{abstract}


\section{Introduction and problem setting}
We consider a multivariate approximation problem in average case setting for special random fields with arbitrary large parametric dimension.

Let $X=\{X(t),\, t\in\R \}$ be a random process defined on some probability space. Here and below $\R$ denotes the set of real numbers. Suppose that the process has zero mean and the following covariance function
\begin{eqnarray*}
	\CovFunc_\s(t,s)=\exp\biggl\{ -\dfrac{(t-s)^2}{2\s^2}\biggr\},\quad t,s\in\R,
\end{eqnarray*}
where $\sigma>0$ is a length scale parameter. The process is usually considered as a random element of the space $L_2(\R, \mu)$, where $\mu$ is the standard Gaussian measure on $\R$. Covariance operator acts as follows
\begin{eqnarray}\label{def_Ksig}
	K_\sigma f(t)=\int\limits_{\R} \CovFunc_{\sigma}(t,s) f(s) \mu(\dd s)=\int\limits_{\R} \CovFunc_\sigma(t,s) f(s) \, \dfrac{e^{-\tfrac{s^2}{2}}}{\sqrt{2\pi}}\, \dd s,\quad t\in\R.
\end{eqnarray}

 We consider $d$-variate version of $X$ with arbitrary large $d\in\N$ (set of positive integers). Namely, we consider a zero-mean random field $Y_d=\bigl\{Y_d(t),\, t\in\R^d\bigr\}$ with the following covariance function
\begin{eqnarray}\label{def_CovFunc_Yd}
	\CovFunc^{Y_d}(t,s)=\prod\limits_{j=1}^{d}\CovFunc_{\s_j}(t_j,s_j)=\exp\biggl\{-\sum\limits_{j=1}^{d} \dfrac{(t_j-s_j)^2}{2\s_j^2}\biggr\},
\end{eqnarray}
where $t=(t_1,\ldots, t_d)$ and $s=(s_1,\ldots, s_d)$ are from $\R^d$. Here $(\s_j)_{j\in\N}$ is a given sequence of length scale parameters, which are generally have different values. If every $\CovFunc_{\s_j}$ corresponds to a zero-mean process $X_j=\{X_j(t),\, t\in\R\}$ (defined on some probability space), $j\in\N$, then  $Y_d$ is called \textit{tensor product of} $X_1$, \ldots, $X_d$ (see \cite{KarNazNik}). Function \eqref{def_CovFunc_Yd} is well known as \textit{Gaussian kernel}, which
is often used in numerical computation and statistical learning (see \cite{Buhmann}, \cite{HastieTibshFried}, \cite{RasWill}, \cite{SchSmola}, \cite{Wendland}). 

For every $d\in\N$ the random field $Y_d$ is considered as random element of the space $L_2(\R^d, \mu_d)$, where $\mu_d$ is the standard Gaussian measure on $\R^d$. So the space is
equipped with the inner product
\begin{eqnarray*}
	\langle f,g \rangle_{2,d}=\int\limits_{\R^d} f(x) g(x)\mu_d(\dd x)=\int\limits_{\R^d} f(x) g(x)\, \dfrac{1}{(2\pi)^{d/2}}\exp\biggl\{-\dfrac{1}{2}\sum\limits_{j=1}^{d} x_j^2  \biggr\} \dd x,
\end{eqnarray*}
and the norm
\begin{eqnarray*}
	\| f \|_{2,d}=\Biggl(\,\,\,\int\limits_{\R^d} f(x)^2\,\mu_d(\dd x)\Biggr)^{1/2}=\Biggl(\,\,\,\int\limits_{\R^d} f(x)^2\, \dfrac{1}{(2\pi)^{d/2}}\exp\biggl\{-\dfrac{1}{2}\sum\limits_{j=1}^{d} x_j^2  \biggr\} \dd x\Biggr)^{1/2},
\end{eqnarray*}
where $x=(x_1,\ldots, x_d)\in\R^d$ in the integrals. The covariance operator $K^{Y_d}$ of $Y_d$ acts as follows  
\begin{eqnarray*}
	K^{Y_d} f(t)=\int\limits_{\R^d} \CovFunc^{Y_d}(t,s) f(s) \mu_d(\dd s)=\int\limits_{\R^d} \CovFunc^{Y_d}(t,s) f(s) \, \dfrac{1}{(2\pi)^{d/2}}\exp\biggl\{-\dfrac{1}{2}\sum\limits_{j=1}^{d} s_j^2  \biggr\} \dd s,
\end{eqnarray*}
where $t=(t_1,\ldots, t_d)$ and $s=(s_1,\ldots, s_d)$ are from $\R^d$.

We consider the \textit{average case approximation complexity} (\textit{approximation complexity} for short) of $Y_d$, $d\in\N$:
\begin{eqnarray}\label{def_nYde}
n^{Y_d}(\e)\colonequals\min\bigl\{n\in\N:\, e^{Y_d}(n)\leqslant \e\, e^{Y_d}(0)  \bigr \},
\end{eqnarray}
where $\e\in(0,1)$ is a given error threshold, and
\begin{eqnarray*}
	e^{Y_d}(n)\colonequals\inf\biggl\{ \Bigl(\Expec\bigl\|Y_d - Y^{( n)}_d\bigr\|_{2,d}^2\Bigr)^{1/2} : Y^{(n)}_d\in \AppClass_n^{Y_d}\biggr\}
\end{eqnarray*}
is the smallest 2-average error among all linear approximations of $Y_d$  having rank $n\in\N$ ($\Expec$ is the expectation). The corresponding classes  of linear algorithms are
\begin{eqnarray*}
	\AppClass_n^{Y_d}\colonequals \Bigl\{\sum_{m=1}^{n} \langle Y_d,\psi_m\rangle_{2,d}\,\psi_m :  \psi_m \in L_2(\R^d,\mu_d)\Bigr \}\,.
\end{eqnarray*}
We will deal with the \textit{normalized error}, i.e. we take into account the quantity:
\begin{eqnarray*}
	e^{Y_d}(0)\colonequals  \bigl(\Expec\|Y_d\|_{2,d}^2\bigr)^{1/2}<\infty,
\end{eqnarray*}
which is the approximation error of $Y_d$ by zero element.

For a given sequence $(\sigma_j)_{j\in\N}$ of length scale parameters in \eqref{def_CovFunc_Yd} the quantity $n^{Y_d}(\e)$ is considered as a function depending on two variables $d\in\N$ and $\e\in(0,1)$. There are a lot of results in this direction concerning the \textit{tractability} (see \cite{NovWoz1}). They provide necessary and sufficient conditions on $(\sigma_j)_{j\in\N}$ to have upper bounds of given forms for the approximation complexity. The results within the described average case setting can be find in the papers \cite{ChenWang}, \cite{FassHickWoz1}, and \cite{Khart2}. The other setting of the worst case was considered  in \cite{FassHickWoz2}, \cite{NovWoz3}, and \cite{SloanWoz}. We will investigate $n^{Y_d}(\e)$ in the different way. Namely, we are interested in the asymptotic behaviour of $n^{Y_d}(\e)$ for arbirarily small fixed $\e$ and $d\to\infty$. We are not aware of any asymptotic results in this way specially for random fields with covariance functions \eqref{def_CovFunc_Yd}. There exist a suitable general methods and results from \cite{Khart1}, but their application requires an additional analysis. So in fact we will do such analysis in this paper. 
 
We will use the following notation. Let $\N_0$ denote the set of non-negative integers. We write $a_n\sim b_n$ if $a_n/b_n\to 1$, $n\to\infty$. The indicator $\id(A)$ equals one if $A$ is true and zero if $A$ is false. For any function $f$ we will denote by $\ContSet(f)$ the set of all its continuity points and by $f^{-1}$ the generalized inverse function $f^{-1}(y)\colonequals\inf\bigl\{x\in\R: f(x)\geqslant y\bigr\}$, where $y$ is from the range of $f$. By \textit{distribution function} $F$ we mean a non-decreasing function $F$ on $\R$ that is right-continuous on $\R$, $\lim\limits_{x\to-\infty} F(x)=0$, and $\lim\limits_{x\to\infty} F(x)=1$. 

\section{Preliminaries}
The quantity $n^{Y_d}(\e)$ can be described in terms of the eigenvalues of the covariance operator $K^{Y_d}$. Let $(\lambda^{Y_d}_m)_{m\in\N}$ denote the sequence of eigenvalues and $(\psi^{Y_d}_m)_{m\in\N}$ the corresponding  sequence of orthonormal eigenvectors of $K^{Y_d}$. The family $(\lambda^{Y_d}_m)_{m\in\N}$ is assumed to be ranked in non-increasing order. We have therefore $K^{Y_d} \psi^{Y_d}_m(t)=\lambda^{Y_d}_m\psi^{Y_d}_m(t)$, $m\in\N$, $t\in\R^d$. We denote by $\Lambda^{Y_d}$ the trace of $K^{Y_d}$, i.e. $\Lambda^{Y_d}\defeq \sum_{m=1}^\infty \lambda^{Y_d}_m$.

It is well known (see \cite{Brown}, \cite{Rit}, \cite{WasWoz}) that for any $n\in\N$ the following $n$-rank random field
\begin{eqnarray}\label{def_Ydn}
\widetilde Y^{(n)}_d(t)\colonequals\sum_{k=1}^n \langle Y_d,\psi^{Y_d}_k\rangle_{2,d}\, \psi^{Y_d}_k(t),\quad t\in\R^d,
\end{eqnarray}
minimizes the 2-average case error.  Hence formula \eqref{def_nYde} is reduced to
\begin{eqnarray*}
	n^{Y_d}(\e)=\min\Bigl\{n\in\N:\, \Expec 
	\bigl\|Y_d-\widetilde{Y}^{(n)}_d\bigr\|_{2,d}^2\leqslant\e^2\,\Expec \|Y_d\|_{2,d}^2 \Bigr\},\quad  d\in\N,\,\, \e\in(0,1).
\end{eqnarray*}
Due to  \eqref{def_Ydn} and $\Expec \langle Y_d,\psi^{Y_d}_m\rangle_{2,d}^2=\langle\psi^{Y_d}_m, K^{Y_d}\psi^{Y_d}_m\rangle_{2,d}=\lambda^{Y_d}_m$, $m\in\N$, we have the needed representation:
\begin{eqnarray*}
n^{Y_d}(\e)=\min\Bigl\{n\in\N:\, \sum_{m=n+1}^\infty \lambda^{Y_d}_m\leqslant\e^2\,\Lambda^{Y_d} \Bigr\},\quad  d\in\N,\,\, \e\in(0,1).
\end{eqnarray*}

We now consider the sequence $(\lambda^{Y_d}_m)_{m\in\N}$. It has the following description. Let $(\lambda_{\sigma,k})_{k\in\N}$ denote the sequence of eigenvalues (ranked in non-increasing order) of the covariance operator $K_\sigma$ defined by  \eqref{def_Ksig}. This sequence is known (see \cite{NovWoz3}, \cite{RasWill}, and \cite{Zhu}):
\begin{eqnarray}\label{def_lambdasigk}
	\lambda_{\sigma,k}=(1-\omega)\,\omega^{k-1}, \quad k\in\N,\quad\text{where}\quad\omega\colonequals \biggr(1+\tfrac{\sigma^2}{2} \Bigl(1+\sqrt{1+\tfrac{4}{\sigma^2}}\,\Bigr)\biggr)^{-1}.
\end{eqnarray}
In particular, we have $\lambda_{\sigma, 1}=1-\omega$ and $\sum_{k\in\N} \lambda_{\sigma, k}=1$.
It is well known (see \cite{LifPapWoz1} and \cite{NovWoz3}) that, due to the tensor product structure \eqref{def_CovFunc_Yd} with given $\sigma_j$, $j\in\N$,   $(\lambda^{Y_d}_m)_{m\in\N}$ is the sequence of numbers
\begin{eqnarray*}
	\prod\limits_{j=1}^{d} \lambda_{\sigma_j,k_j}=\prod\limits_{j=1}^{d} (1-\omega_j)\,\omega_j^{k_j-1},\quad k_1,k_2,\ldots, k_d \in\N,
\end{eqnarray*}
ranked in non-increasing order (see \cite{NovWoz1}). Here, according to \eqref{def_lambdasigk}, we set
\begin{eqnarray}\label{def_omegaj}
	\omega_j\colonequals \biggl(1+\tfrac{\sigma_j^2}{2} \Bigl(1+\sqrt{1+\tfrac{4}{\sigma_j^2}}\,\Bigr)\biggr)^{-1},\quad j\in\N.
\end{eqnarray}
Observe that
\begin{eqnarray*}
	\Lambda^{Y_d}=\prod\limits_{j=1}^{d} \sum\limits_{k\in\N}\lambda_{\sigma_j,k}=1,\quad d\in\N.
\end{eqnarray*}

Thus each of the sequences $(\sigma_j)_{j\in\N}$ and $(\omega_j)_{j\in\N}$ fully determines $(\lambda^{Y_d}_m)_{m\in\N}$ and hence $n^{Y_d}(\e)$ for any  $d\in\N$ and $\e\in(0,1)$.

\section{General results}
Before proceeding to the asymptotic analysis of the quantity $n^{Y_d}(\e)$, we find criteria of its boundedness and unboundedness on $d$ for any fixed $\e\in(0,1)$. The following propositions show that for any fixed $\e\in(0,1)$ either the quantity $n^{Y_d}(\e)$ is a bounded function on $d\in\N$ or it tends to infinity as $d\to\infty$.

\begin{Proposition}	\label{pr_boundedness}
	The following conditions are equivalent:
	\begin{itemize}
		\item[$(i)$] $\sup_{d \in \mathbb{N}} n^{Y_d}(\e) < \infty$ for every $\e\in(0, 1)$;
		\item[$(ii)$] $\sum_{j = 1}^\infty \omega_j< \infty$;
		\item[$(iii)$] $\sum_{j = 1}^\infty \sigma_j^{-2} < \infty$.
	\end{itemize}
\end{Proposition}
\textbf{Proof of Proposition \ref{pr_boundedness}.}\quad By Proposition 5 from \cite{Khart1}, the relation $\sup_{d \in\N} n^{Y_d}(\e)<\infty$, $\e\in(0,1)$, is equivalent to convergence of the following series
\begin{eqnarray}\label{def_series_boundedness}
	\sum_{j = 1}^\infty\sum_{k=2}^\infty \dfrac{\lambda_{\sigma_j,k}}{\lambda_{\sigma_j,1}}=\sum_{j = 1}^\infty\sum_{k=2}^\infty\omega_j^{k-1}=\sum_{j = 1}^\infty\dfrac{\omega_j}{1 - \omega_j}=\sum_{j = 1}^\infty  \dfrac{2}{\sigma_j^2+\sqrt{\sigma_j^4+4\sigma_j^{2}}}.
\end{eqnarray}
Since $\omega_j\in(0,1)$, the convergence of  $\sum_{j = 1}^\infty\tfrac{\omega_j}{1 - \omega_j}$ implies the convergence of $\sum_{j = 1}^\infty\omega_j$. Next, if $\sum_{j = 1}^\infty\omega_j<\infty$, then $\omega_j\to 0$ and hence $\tfrac{\omega_j}{1 - \omega_j}\sim \omega_j$, $j\to\infty$. So we have $\sum_{j = 1}^\infty \tfrac{\omega_j}{1 - \omega_j}<\infty$.

It is easily seen that the convergence of $\sum_{j = 1}^\infty \sigma_j^{-2}$ implies the convergence of \eqref{def_series_boundedness}. Next, if \eqref{def_series_boundedness} converges, then $\sigma_j\to \infty$ and $\sigma_j^2 + \sqrt{\sigma_j^4+4\sigma_j^2} \sim 2 \sigma_j^2$, $j\to\infty$. Hence we get the convergence of $\sum_{j = 1}^\infty \sigma_j^{-2}$. \quad $\Box$\\

\begin{Proposition}	\label{pr_unboundedness}
	The following conditions are equivalent:
	\begin{itemize}
		\item[$(i)$] $\lim_{d \to\infty} n^{Y_d}(\e) = \infty$ for every $\e\in(0, 1)$;
		\item[$(ii)$] $\sum_{j = 1}^\infty \omega_j=\infty$;
		\item[$(iii)$] $\sum_{j = 1}^\infty \sigma_j^{-2} = \infty$.
	\end{itemize}
\end{Proposition}
\textbf{Proof of Proposition \ref{pr_unboundedness}.}\quad According to Proposition 4 and 5 from \cite{Khart1}, for any $\e\in(0,1)$ either $\sup_{d \in\N} n^{Y_d}(\e)<\infty$ or  $\lim_{d \to\infty} n^{Y_d}(\e) = \infty$. Hence the latter is equivalent to divergence of each of the series $\sum_{j = 1}^\infty \omega_j$ and $\sum_{j = 1}^\infty \sigma_j^{-2}$ by Proposition \ref{pr_boundedness}.\quad $\Box$\\

It is known (see \cite{Khart1}) that for wide class of tensor product random fields the quantity $n^{Y_d}(\e)$ has the logarithmic asymptotics of the form \eqref{th_nYde_Limitdistr_cond} below.  Our next theorem shows that  the function $q$ can be  only a special quantile of self-decomposable distribution function in such asymptotics (see \cite{GnedKolm}, \cite{Petrov}, or \cite{Khart1}, Appendix). Recall that self-decomposable distribution functions are completely described by the triplet $(c,v, L)$ from spectral representation of their characteristic functions, where $c\in\R$ is a shift parameter, $v>0$ is the Gaussian component, $L$ is the L\'evy spectral function (see \cite{Khart1}, Appendix).

\begin{Theorem}\label{th_nYde_Limitdistr}
	Let $(a_d)_{d\in\N}$ be a sequence, $(b_d)_{d\in\N}$ be a positive sequence such that $b_d\to\infty$, $d\to\infty$. Let  a non-increasing function $q\colon (0,1)\to\R$ and a distribution function $G$ satisfy the equation $q(\e)=G^{-1}(1-\e^2)$ for all $\e\in\ContSet(q)$. Suppose that the following asymptotics holds
	\begin{eqnarray}\label{th_nYde_Limitdistr_cond}
		\forall \e\in\ContSet(q)\quad\ln n^{Y_d}(\e)=a_d+q(\e) b_d+o(b_d),\quad d\to\infty.
	\end{eqnarray}
    Then $G$ is self-decomposable with zero L\'evy spectral function on $(-\infty,0)$. 
\end{Theorem}
\textbf{Proof of Theorem \ref{th_nYde_Limitdistr}.}\quad Due to Theorem 1 from \cite{Khart1}, the condition \eqref{th_nYde_Limitdistr_cond} is equivalent to the convergence
\begin{eqnarray}\label{conc_Gd_conv}
	\lim_{d\to\infty} G_d(x)= G(x),\quad x\in \ContSet(G),
\end{eqnarray}
where
\begin{eqnarray}\label{def_Gd}
	G_d(x)\defeq \sum\limits_{m\in\N} \lambda_m^{Y_d}\,\id\bigl(\lambda_m^{Y_d}\geqslant e^{-a_d-b_dx}\bigr),\quad x\in\R,\quad d\in\N.  
\end{eqnarray}
It was shown in \cite{Khart1} that $G_d$, $d\in\N$, can be considered as the following distribution functions
\begin{eqnarray*}
	G_d(x)=\Probab \Biggl(\dfrac{\sum_{j=1}^d U_j - a_d}{b_d} \leqslant x\Biggr), \quad x\in\R,\quad d\in\N, 
\end{eqnarray*} 
where  $U_j$, $j\in\N$, are independent random variables on some probability space with the measure $\Probab$. Here $U_j$, $j\in\N$, have the following distribution
\begin{eqnarray*}
	\Probab\bigl(U_j=|\ln \lambda_{\sigma_j, k}| \bigr)=\lambda_{\sigma_j, k},\quad k\in\N,\quad j\in\N.
\end{eqnarray*}

Now we center these variables in the following way: $\hat U_j\defeq U_j-|\ln\lambda_{\sigma_j,1}|$, $j\in\N$. So we have
\begin{eqnarray*}
	\Probab\bigl(\hat U_j=|\ln \lambda_{\sigma_j, k}|-|\ln \lambda_{\sigma_j,1}| \bigr)=\lambda_{\sigma_j, k},\quad k\in\N,\quad j\in\N,
\end{eqnarray*}
i.e.
\begin{eqnarray*}
	\Probab\bigl(\hat U_j=k|\ln \omega_j|\bigr)=(1-\omega_j)\omega_j^k,\quad k\in\N_0,\quad j\in\N.
\end{eqnarray*}
Here $\hat U_j\geqslant 0$, $j\in\N$. Next, we set
\begin{eqnarray*}
	\hat a_d\defeq a_d-\sum_{j=1}^d|\ln \lambda_{\sigma_j,1}| = a_d -\sum_{j=1}^d |\ln(1-\omega_j)|,\quad d\in\N.
\end{eqnarray*}
Then
\begin{eqnarray}\label{conc_Gdx}
	G_d(x)=\Probab \Biggl(\dfrac{\sum_{j=1}^d \hat U_j - \hat a_d}{b_d} \leqslant x\Biggr), \quad x\in\R,\quad d\in\N, 
\end{eqnarray}

For any $d\in\N$, $j\in\{1,\ldots, d\}$ and $x>0$ we consider the following distribution tails:
\begin{eqnarray*}
	\Probab(|\hat U_j|>x b_d)=\Probab(\hat U_j>x b_d)=\sum\limits_{\substack{k\in\N_0:\\k|\ln \omega_j|>x b_d}} (1-\omega_j)\omega_j^k = \omega_j^{k_{j,d}(x)},
\end{eqnarray*}
where $k_{j,d}(x)\defeq\min\{k\in\N_0:  k|\ln \omega_j|>x b_d\}$. Since $k_{j,d}(x)|\ln \omega_j|>x b_d$, we have
\begin{eqnarray*}
	\Probab(|\hat U_j|>x b_d)<e^{-x b_d},\quad d\in\N,\quad j\in\{1,\ldots, d\},\quad x>0.
\end{eqnarray*}
Due to $b_d\to\infty$, $d\to\infty$,  we obtain
\begin{eqnarray}\label{cond_uni_neglig_hUj}
	\max_{j\in\{1,\ldots, d\}}\Probab(|\hat U_j|>x b_d)\to 0,\quad d\to \infty.
\end{eqnarray}
This is the condition of uniform negligibility of $\hat U_j/b_d$ in the sums   $\bigl(\sum_{j = 1}^n \hat U_j -\hat a_d\bigr)/b_d$. It is known, that under this condition and \eqref{conc_Gdx} the limit distribution in \eqref{conc_Gd_conv} is self-decomposable (see \cite{Petrov}, p. 101, or \cite{Khart1}, Theorem 10). Let $L$ denote the L\'evy spectral function of $G$. Due to the non-negativity of $\hat U_j$, we have $L(x)=0$, $x<0$ (see \cite{GnedKolm}, p. 124, or \cite{Khart1}, Theorem 11).\quad $\Box$\\

The next theorem provides a criterion for the asymptotics \eqref{th_nYde_Limitdistr_cond}, where $q$ is a quantile of a given self-decomposable distribution function. This is the main result of the paper.

\begin{Theorem}	\label{th_nYde_LogAsymp}
	Let $(a_d)_{d \in \mathbb{N}}$ be a sequence, $(b_d)_{d \in \mathbb{N}}$ be a positive sequence such that $b_d \to +\infty$, $d \to \infty$. Let  $G$
	be a self-decomposable distribution function with triplet $(\cen, \var, L)$ such that $L(x) = 0$, $x < 0$.  Let  a non-increasing function $q\colon (0,1)\to\R$ satysfy the equation $q(\e)=G^{-1}(1-\e^2)$ for all $\e\in(0,1)$. For the asymptotics
	\begin{eqnarray}\label{th_nYde_LogAsymp_cond}
	\forall \e\in (0,1)\quad\ln n^{Y_d}(\e)=a_d+q(\e) b_d+o(b_d),\quad d\to\infty,
    \end{eqnarray}
    the following ensemble of conditions is necessary and sufficient:
\begin{eqnarray*}
	\mathrm{(A)}&&\lim_{d \to \infty} \sum_{\substack{j = 1,\ldots, d\\ |\ln \omega_j| > \tau b_d}} \omega_j = -L(\tau),\quad \tau>0;\\
	\mathrm{(B)}&&\lim_{d \to \infty}  \frac{1}{b_d} \biggl(\sum_{\substack{j = 1,\ldots, d\\ |\ln \omega_j| \leqslant \tau b_d}} \frac{ |\ln \omega_j|\omega_j}{1 - \omega_j} - \hat{a}_d \biggr) =
	\cen + \gamma_\tau,\quad \tau > 0; \\
	\mathrm{(C)}&&\lim_{\tau \to 0} \varliminf_{d \to \infty} \frac{1}{b_d^2} \sum_{\substack{j = 1,\ldots, d\\ |\ln \omega_j| \leqslant \tau b_d}} \frac{ |\ln \omega_j|^2\omega_j}{(1 - \omega_j)^2} =
	\lim_{\tau \to 0} \varlimsup_{d \to \infty} \frac{1}{b_d^2} \sum_{\substack{j = 1,\ldots, d\\ |\ln \omega_j| \leqslant \tau b_d}} \frac{ |\ln \omega_j|^2\omega_j}{(1 - \omega_j)^2}  = \var; 
\end{eqnarray*}
where 
\begin{eqnarray}\label{def_gammatau}
	&&\gamma_\tau\defeq \int\displaylimits_0^\tau \frac{y^3\dd L(y)}{1 + y^2} - \int\displaylimits_\tau^{+\infty}\frac{y\dd L(y)}{1 + y^2},\quad \tau>0,\\
	\label{def_hatad}&&\hat{a}_d \defeq a_d - \sum_{j = 1}^d|\ln(1 - \omega_j)|,\quad d\in\N.
\end{eqnarray}
\end{Theorem}

The proof of this theorem is essentially based on the following lemma.

\begin{Lemma}\label{lm_sum_omegaj} 
For any $x>0$ and $d\in\N$ the following identities hold:
\begin{eqnarray*}
&&\sum_{j = 1}^d \sum_{\substack{k \in\N:\\ k |\ln \omega_j| > x}} (1 - \omega_j) \omega_j^k = \sum_{\substack{j = 1,\ldots,d:\\ |\ln \omega_j| > x}} \omega_j  + R_0(d,x),\\
&&\sum_{j = 1}^d \sum_{\substack{k \in\N:\\ k |\ln \omega_j| \leqslant x}} k |\ln \omega_j| (1 - \omega_j) \omega_j^k =\sum_{\substack{j = 1,\ldots,d:\\ |\ln \omega_j| \leqslant x}}\dfrac{ |\ln \omega_j|\,\omega_j}{1 - \omega_j}-R_1(d,x),\\
&&\sum_{j = 1}^d \Biggl[\sum_{\substack{k \in\N:\\ k |\ln \omega_j| \leqslant x}} k^2 |\ln \omega_j|^2 (1 - \omega_j) \omega_j^k- \biggl(\sum_{\substack{k \in\N:\\ k |\ln \omega_j| \leqslant x}} k |\ln \omega_j| (1 - \omega_j) \omega_j^k\biggr)^2\Biggr]\\
&&\qquad\qquad\qquad\qquad\qquad\qquad\qquad\qquad\qquad\qquad\qquad=\sum_{\substack{j = 1,\ldots,d:\\ |\ln \omega_j| \leqslant x}}\dfrac{ |\ln \omega_j|^2\omega_j}{(1 - \omega_j)^2}-R_2(d,x),
\end{eqnarray*}
where 
\begin{eqnarray*}
	k_j(x) = \min\bigl\{k\in\N : k\geqslant 2,\, k |\ln \omega_j| > x \bigr\},
\end{eqnarray*}
and
\begin{eqnarray*}
	R_{0}(d,x)&\defeq& \sum_{\substack{j = 1,\ldots,d:\\ |\ln \omega_j| \leqslant x}} \omega_j^{k_{j}(x)},\\
	R_1(d,x)&\defeq&\sum_{\substack{j = 1,\ldots,d:\\ |\ln \omega_j| \leqslant x}} \dfrac{|\ln \omega_j|\,\omega_j^{k_{j}(x)}}{1 - \omega_j} \bigl(k_{j}(x)(1 - \omega_j) + \omega_j\bigr),\\
	R_2(d,x)&\defeq& \sum_{\substack{j = 1,\ldots,d:\\ |\ln \omega_j| \leqslant x}}\dfrac{|\ln \omega_j|^2\omega_j^{k_j(x)}}{(1-\omega_j)^2}\biggl(k_j(x)^2(1-\omega_j)^2\bigl(1+\omega_j^{k_j(x)}\bigr)\\
	&&{}\quad+2k_j(x)(1-\omega_j)\omega_j^{k_j(x)+1}+(1-\omega_j)\omega_j +\omega_j^{k_j(x)+2}\biggr).
\end{eqnarray*}
\end{Lemma}
\textbf{Proof of Lemma \ref{lm_sum_omegaj}.}\quad We fix $x>0$ and $d\in\N$. 

1. Let us prove the first identity. Observe that
\begin{eqnarray*}
	\sum_{j = 1}^d \sum_{\substack{k \in\N:\\ k |\ln \omega_j| > x}} (1 - \omega_j) \omega_j^k=\sum_{\substack{j = 1,\ldots,d:\\ |\ln \omega_j| > x}}\, \sum_{\substack{k \in\N:\\ k |\ln \omega_j| > x}} (1 - \omega_j) \omega_j^k+\sum_{\substack{j = 1,\ldots,d:\\ |\ln \omega_j| \leqslant x}}\, \sum_{\substack{k \in\N:\\ k |\ln \omega_j| > x}} (1 - \omega_j) \omega_j^k.
\end{eqnarray*}
Here we have
\begin{eqnarray*}
	\sum_{\substack{j = 1,\ldots,d:\\ |\ln \omega_j| > x}}\, \sum_{\substack{k \in\N:\\ k |\ln \omega_j| > x}} (1 - \omega_j) \omega_j^k=\sum_{\substack{j = 1,\ldots,d:\\ |\ln \omega_j| > x}} \sum_{k \in\N} (1 - \omega_j) \omega_j^k=\sum_{\substack{j = 1,\ldots,d:\\ |\ln \omega_j| > x}} \omega_j,
\end{eqnarray*}
and
\begin{eqnarray*}
	\sum_{\substack{j = 1,\ldots,d:\\ |\ln \omega_j| \leqslant x}} \sum_{\substack{k \in\N:\\ k |\ln \omega_j| > x}} (1 - \omega_j) \omega_j^k&=&\sum_{\substack{j = 1,\ldots,d:\\ |\ln \omega_j| \leqslant x}} \sum_{\substack{k \in\N:\, k\geqslant 2,\\ k |\ln \omega_j| > x}} (1 - \omega_j) \omega_j^k\\
	&=&\sum_{\substack{j = 1,\ldots,d:\\ |\ln \omega_j| \leqslant x}} \sum_{k=k_j(x)}^\infty (1 - \omega_j) \omega_j^k=\sum_{\substack{j = 1,\ldots,d:\\ |\ln \omega_j| \leqslant x}}\omega_j^{k_j(x)}.
\end{eqnarray*}

2. Let us prove the second identity.   It is obvious that if $|\ln \omega_j|>x$ then there are no any $k\in\N$ such that $k|\ln \omega_j|\leqslant x$. Therefore
\begin{eqnarray*}
	\sum_{j = 1}^d \sum_{\substack{k \in\N:\\ k |\ln \omega_j| \leqslant x}} k |\ln \omega_j| (1 - \omega_j) \omega_j^k&=&\sum_{\substack{j = 1,\ldots,d:\\ |\ln \omega_j| \leqslant x}} \sum_{\substack{k \in\N:\\ k |\ln \omega_j| \leqslant x}} k |\ln \omega_j| (1 - \omega_j) \omega_j^k\\
	&=&\sum_{\substack{j = 1,\ldots,d:\\ |\ln \omega_j| \leqslant x}} |\ln \omega_j| (1 - \omega_j) \sum_{\substack{k \in\N:\\ k |\ln \omega_j| \leqslant x}} k \omega_j^k.
\end{eqnarray*}
Here
\begin{eqnarray*}
	\sum_{\substack{k \in\N:\\ k |\ln \omega_j| \leqslant x}} k \omega_j^k&=&\sum_{k=1}^\infty k  \omega_j^k-\sum_{k=k_j(x)}^\infty k  \omega_j^k.
\end{eqnarray*}

It is well known that
\begin{eqnarray*}
	\sum_{k=1}^\infty k  \omega_j^k= \omega_j\sum_{k=1}^\infty k  \omega_j^{k-1}=\dfrac{\omega_j}{(1-\omega_j)^2}.
\end{eqnarray*}
Next, using this fact we get
\begin{eqnarray*}
	\sum_{k=k_j(x)}^\infty k \omega_j^k&=&\omega_j^{k_j(x)}\sum_{k=k_j(x)}^\infty k \omega_j^{k-k_j(x)}\\
	&=& \omega_j^{k_j(x)}	\sum_{k=k_j(x)+1}^\infty (k-k_j(x)) \omega_j^{k-k_j(x)}+k_j(x)\omega_j^{k_j(x)}\sum_{k=k_j(x)}^\infty \omega_j^{k-k_j(x)}\\
	&=&\dfrac{\omega_j^{k_j(x)+1}}{(1-\omega_j)^2}+\dfrac{k_j(x)\omega_j^{k_j(x)}}{1-\omega_j}.
\end{eqnarray*}
Then
\begin{eqnarray*}
	\sum_{j = 1}^d \sum_{\substack{k \in\N:\\ k |\ln \omega_j| \leqslant x}} k |\ln \omega_j| (1 - \omega_j) \omega_j^k
	&=&\sum_{\substack{j = 1,\ldots,d:\\ |\ln \omega_j| \leqslant x}} |\ln \omega_j| (1 - \omega_j)\Biggl(\dfrac{\omega_j}{(1-\omega_j)^2}-\dfrac{\omega_j^{k_j(x)+1}}{(1-\omega_j)^2}-\dfrac{k_j(x)\omega_j^{k_j(x)}}{1-\omega_j}\Biggr)\\
	&=&\sum_{\substack{j = 1,\ldots,d:\\ |\ln \omega_j| \leqslant x}} \dfrac{|\ln \omega_j|\omega_j}{1-\omega_j}-\sum_{\substack{j = 1,\ldots,d:\\ |\ln \omega_j| \leqslant x}} \Biggl( \dfrac{|\ln \omega_j|\omega_j^{k_j(x)+1}}{1-\omega_j}+k_j(x)|\ln \omega_j|\omega_j^{k_j(x)}\Biggr)\\
	&=&\sum_{\substack{j = 1,\ldots,d:\\ |\ln \omega_j| \leqslant x}} \dfrac{|\ln \omega_j|\omega_j}{1-\omega_j}-\sum_{\substack{j = 1,\ldots,d:\\ |\ln \omega_j| \leqslant x}}\dfrac{|\ln \omega_j|\omega_j^{k_j(x)}}{1-\omega_j} \bigl(\omega_j +k_j(x)(1-\omega_j)\bigr).
\end{eqnarray*}

3. We now prove the third identity. Let us consider the sum
\begin{eqnarray*}
	S_1&\defeq& \sum_{j = 1}^d \sum_{\substack{k \in\N:\\ k |\ln \omega_j| \leqslant x}} k^2 |\ln \omega_j|^2 (1 - \omega_j) \omega_j^k
	=\sum_{\substack{j = 1,\ldots,d:\\ |\ln \omega_j| \leqslant x}}\sum_{\substack{k \in\N:\\ k |\ln \omega_j| \leqslant x}} k^2 |\ln \omega_j|^2 (1 - \omega_j) \omega_j^k\\
	&&{}=\sum_{\substack{j = 1,\ldots,d:\\ |\ln \omega_j| \leqslant x}}|\ln \omega_j|^2(1 - \omega_j)  \sum_{\substack{k \in\N:\\ k |\ln \omega_j| \leqslant x}} k^2 \omega_j^k.
\end{eqnarray*}
Here
\begin{eqnarray*}
	\sum_{\substack{k \in\N:\\ k |\ln \omega_j| \leqslant x}} k^2 \omega_j^k=\sum_{k =1}^\infty k^2  \omega_j^k-\sum_{k=k_j(x)}^\infty k^2  \omega_j^k.
\end{eqnarray*}

It is easily seen that
\begin{eqnarray*}
	\sum_{k=1}^\infty k^2  \omega_j^k&=&\sum_{k=1}^\infty k(k-1)  \omega_j^{k}+\sum_{k =1}^\infty k  \omega_j^{k}
	=\omega_j^2\sum_{k=2}^\infty k(k-1)  \omega_j^{k-2}+\sum_{k=1}^\infty k  \omega_j^{k}\\
	&&{}=\dfrac{2\omega_j^2}{(1-\omega_j)^3}+\dfrac{\omega_j}{(1-\omega_j)^2}
	= \dfrac{\omega_j+\omega_j^2}{(1-\omega_j)^3}.
\end{eqnarray*}
Next, we write
\begin{eqnarray*}
	\sum_{k=k_j(x)}^\infty k^2 \omega_j^k
	&=& \sum_{k=k_j(x)}^\infty (k-k_j(x))^2 \omega_j^{k}+2k_j(x)\sum_{k=k_j(x)}^\infty k\omega_j^{k}-k_j(x)^2\sum_{k=k_j(x)}^\infty \omega_j^{k}\\
	&=&\omega_j^{k_j(x)}\sum_{k=k_j(x)+1}^\infty (k-k_j(x))^2 \omega_j^{k-k_j(x)}+2k_j(x)\sum_{k=k_j(x)}^\infty k\omega_j^{k}-k_j(x)^2\sum_{k=k_j(x)}^\infty \omega_j^{k}.
\end{eqnarray*}
Due to above remarks we get
\begin{eqnarray*}
	\sum_{k=k_j(x)}^\infty k^2 \omega_j^k&=&\omega_j^{k_j(x)}\Biggl(\dfrac{2\omega_j^2}{(1-\omega_j)^3}+\dfrac{\omega_j}{(1-\omega_j)^2}\Biggr)+2k_j(x)\Biggl( \dfrac{\omega_j^{k_j(x)+1}}{(1-\omega_j)^2}+\dfrac{k_j(x)\omega_j^{k_j(x)}}{1-\omega_j}\Biggr)\\
	&&{}-k_j(x)^2\cdot \dfrac{\omega_j^{k_j(x)}}{1-\omega_j}\\
	&=&\dfrac{2\omega_j^{k_j(x)+2}}{(1-\omega_j)^3}+\dfrac{(2k_j(x)+1)\omega_j^{k_j(x)+1}}{(1-\omega_j)^2}+\dfrac{k_j(x)^2\omega_j^{k_j(x)}}{1-\omega_j}.
\end{eqnarray*}
Thus we obtain
\begin{eqnarray*}
	S_1&=&\sum_{\substack{j = 1,\ldots,d:\\ |\ln \omega_j| \leqslant x}}|\ln \omega_j|^2(1 - \omega_j) \Biggl(\dfrac{\omega_j+\omega_j^2}{(1-\omega_j)^3}- \dfrac{2\omega_j^{k_j(x)+2}}{(1-\omega_j)^3}-\dfrac{(2k_j(x)+1)\omega_j^{k_j(x)+1}}{(1-\omega_j)^2}-\dfrac{k_j(x)^2\omega_j^{k_j(x)}}{1-\omega_j}\Biggr)\\
	&=&\sum_{\substack{j = 1,\ldots,d:\\ |\ln \omega_j| \leqslant x}}|\ln \omega_j|^2\Biggl(\dfrac{\omega_j+\omega_j^2}{(1-\omega_j)^2}- \dfrac{2\omega_j^{k_j(x)+2}}{(1-\omega_j)^2}-\dfrac{(2k_j(x)+1)\omega_j^{k_j(x)+1}}{1-\omega_j}-k_j(x)^2\omega_j^{k_j(x)}\Biggr).
\end{eqnarray*}

We now consider the sums
\begin{eqnarray*}
	S_2&\defeq&\sum_{j = 1}^d \biggl(\sum_{\substack{k \in\N:\\ k |\ln \omega_j| \leqslant x}} k |\ln \omega_j| (1 - \omega_j) \omega_j^k\biggr)^2\\
	&=&\sum_{\substack{j = 1,\ldots,d:\\ |\ln \omega_j| \leqslant x}}\biggl(\sum_{\substack{k \in\N:\\ k |\ln \omega_j| \leqslant x}} k |\ln \omega_j| (1 - \omega_j) \omega_j^k\biggr)^2\\
	&=&\sum_{\substack{j = 1,\ldots,d:\\ |\ln \omega_j| \leqslant x}}|\ln \omega_j|^2(1 - \omega_j)^2 \biggl(\sum_{k=1}^\infty k  \omega_j^k-\sum_{k=k_j(x)}^\infty k  \omega_j^k \biggr)^2.
\end{eqnarray*}
By the above remarks, we have
\begin{eqnarray*}
	\biggl(\sum_{k=1}^\infty k  \omega_j^k-\sum_{k=k_j(x)}^\infty k  \omega_j^k\biggr)^2&=&\Biggl(\dfrac{\omega_j}{(1-\omega_j)^2}-\dfrac{\omega_j^{k_j(x)+1}}{(1-\omega_j)^2}-\dfrac{k_j(x)\omega_j^{k_j(x)}}{1-\omega_j}\Biggr)^2\\
	&=&\dfrac{\omega_j^2}{(1-\omega_j)^4}+\dfrac{\omega_j^{2k_j(x)+2}}{(1-\omega_j)^4}+\dfrac{k_j(x)^2\omega_j^{2k_j(x)}}{(1-\omega_j)^2}-\dfrac{2\omega_j^{k_j(x)+2}}{(1-\omega_j)^4}\\
	&&{}-\dfrac{2k_j(x)\omega_j^{k_j(x)+1}}{(1-\omega_j)^3}+\dfrac{2k_j(x)\omega_j^{2k_j(x)+1}}{(1-\omega_j)^3}.
\end{eqnarray*}
Hence we get
\begin{eqnarray*}
	S_2&=&\sum_{\substack{j = 1,\ldots,d:\\ |\ln \omega_j| \leqslant x}}|\ln \omega_j|^2\Biggl(\dfrac{\omega_j^2}{(1-\omega_j)^2}+\dfrac{\omega_j^{2k_j(x)+2}}{(1-\omega_j)^2}+k_j(x)^2\omega_j^{2k_j(x)}-\dfrac{2\omega_j^{k_j(x)+2}}{(1-\omega_j)^2}\\
	&&{}-\dfrac{2k_j(x)\omega_j^{k_j(x)+1}}{1-\omega_j}+\dfrac{2k_j(x)\omega_j^{2k_j(x)+1}}{1-\omega_j} \Biggr).
\end{eqnarray*}

Thus we have
\begin{eqnarray*}
	S_1-S_2&=&\sum_{\substack{j = 1,\ldots,d:\\ |\ln \omega_j| \leqslant x}}|\ln \omega_j|^2\Biggl(\dfrac{\omega_j}{(1-\omega_j)^2}- \dfrac{\omega_j^{k_j(x)+1}}{1-\omega_j}-k_j(x)^2\omega_j^{k_j(x)}\\
	&&{}-\dfrac{\omega_j^{2k_j(x)+2}}{(1-\omega_j)^2}-k_j(x)^2\omega_j^{2k_j(x)}-\dfrac{2k_j(x)\omega_j^{2k_j(x)+1}}{1-\omega_j} \Biggr)\\
	&=&\sum_{\substack{j = 1,\ldots,d:\\ |\ln \omega_j| \leqslant x}}\dfrac{|\ln \omega_j|^2\omega_j}{(1-\omega_j)^2}-\sum_{\substack{j = 1,\ldots,d:\\ |\ln \omega_j| \leqslant x}}\dfrac{|\ln \omega_j|^2\omega_j^{k_j(x)}}{(1-\omega_j)^2}\Biggl((1-\omega_j)\omega_j+k_j(x)^2(1-\omega_j)^2\\
	&&{}+\omega_j^{k_j(x)+2}+k_j(x)^2(1-\omega_j)^2\omega_j^{k_j(x)}+2k_j(x)(1-\omega_j)\omega_j^{k_j(x)+1}\Biggr)\\
	&=&\sum_{\substack{j = 1,\ldots,d:\\ |\ln \omega_j| \leqslant x}}\dfrac{|\ln \omega_j|^2\omega_j}{(1-\omega_j)^2}-\sum_{\substack{j = 1,\ldots,d:\\ |\ln \omega_j| \leqslant x}}\dfrac{|\ln \omega_j|^2\omega_j^{k_j(x)}}{(1-\omega_j)^2}\biggl(k_j(x)^2(1-\omega_j)^2\bigl(1+\omega_j^{k_j(x)}\bigr)\\
	&&{}\quad+2k_j(x)(1-\omega_j)\omega_j^{k_j(x)+1}+(1-\omega_j)\omega_j +\omega_j^{k_j(x)+2}\biggr).
\end{eqnarray*}
So we obtain the third identity. \quad $\Box$\\

We  now turn to proof of the theorem.\\

\noindent
\textbf{Proof of Theorem \ref{th_nYde_LogAsymp}.}\quad We first recall that $\ContSet(q)=(0,1)$, which follows from properties of self-decomposable distribution functions (see \cite{Khart1}, Remarks 2 and 8). Thus, by Theorem~1 from \cite{Khart1}, the condition \eqref{th_nYde_LogAsymp_cond} is equivalent to the convergence \eqref{conc_Gd_conv}, where $G$ is the given self-decomposable law with the triplet $(c,v,L)$ and $G_d$, $d\in\N$, are defined by \eqref{def_Gd} for given $a_d$, $b_d$, and $Y_d$, $d\in\N$. According to \eqref{conc_Gdx}, $G_d$, $d\in\N$, are distribution functions of centered and normalized sums of non-negative independent random variables $\hat{U}_j$, $j\in\N$, satisfying \eqref{cond_uni_neglig_hUj}. Hence for the convergence \eqref{conc_Gd_conv} the conditions ($\mathrm{A_1}$), ($\mathrm{B}$), ($\mathrm{C}$) of Theorem 11 from \cite{Khart1} (Appendix) are necessary and sufficient (with $Y_j\defeq \hat{U}_j$, $j\in\N$, $A_d\defeq \hat{a}_d$, $B_d\defeq b_d$, $\gamma\defeq c$, and $\s^2\defeq v$). Thus \eqref{th_nYde_LogAsymp_cond} is equivalent to the following ensemble of conditions:
\begin{eqnarray*}
	\lim_{d\to\infty}\sum_{j=1}^{d} \Probab(\hat{U}_j>\tau b_d)=-L(\tau)\quad\text{for all}\quad \tau>0;
\end{eqnarray*}
\begin{eqnarray*}
	\lim_{d\to\infty}\dfrac{1}{b_d} \biggl(\sum_{j=1}^{d} \Expec \Bigl[\hat{U}_j\id(|\hat{U}_j|\leqslant \tau b_d)\Bigr]-\hat{a}_d\biggr)=c+\gamma_\tau \quad\text{for all}\quad \tau>0;
\end{eqnarray*}
\begin{eqnarray*}	
	\lim_{\tau\to0}\varliminf_{d\to\infty} \dfrac{1}{b_d^2}\sum_{j=1}^{d}\Var \Bigl[\hat{U}_j\id(|\hat{U}_j|\leqslant \tau b_d)\Bigr]=\lim_{\tau\to0}\varlimsup_{d\to\infty} \dfrac{1}{b_d^2}\sum_{j=1}^{d}\Var \Bigl[\hat{U}_j\id(|\hat{U}_j|\leqslant \tau b_d)\Bigr]=v.
\end{eqnarray*}
Here $\gamma_\tau$ is defined by \eqref{def_gammatau}.  

Let us write the sums in these conditions in terms of $\omega_j$. First, note that 
\begin{eqnarray*}
	\Probab(\hat{U}_j>\tau b_d)&=& \sum_{\substack{k \in\N:\\ k |\ln \omega_j| > \tau b_d}} (1 - \omega_j) \omega_j^k,\\
	\Expec \Bigl[\hat{U}_j\id(|\hat{U}_j|\leqslant \tau b_d)\Bigr]&=&\sum_{\substack{k \in\N:\\ k |\ln \omega_j| \leqslant \tau b_d}} k |\ln \omega_j| (1 - \omega_j) \omega_j^k,\\
	\Var \Bigl[\hat{U}_j\id(|\hat{U}_j|\leqslant \tau b_d)\Bigr]&=&\Expec \Bigl[\hat{U}_j^2\id(|\hat{U}_j|\leqslant \tau b_d)\Bigr]-\biggl(\Expec \Bigl[\hat{U}_j\id(|\hat{U}_j|\leqslant \tau b_d)\Bigr]\biggr)^2\\
	&=&\sum_{\substack{k \in\N:\\ k |\ln \omega_j| \leqslant \tau b_d}} k^2 |\ln \omega_j|^2 (1 - \omega_j) \omega_j^k- \biggl(\sum_{\substack{k \in\N:\\ k |\ln \omega_j| \leqslant \tau b_d}} k |\ln \omega_j| (1 - \omega_j) \omega_j^k\biggr)^2.
\end{eqnarray*}
Next, applying Lemma \ref{lm_sum_omegaj} we have
\begin{eqnarray*}
	\sum_{j=1}^{d} \Probab(\hat{U}_j>\tau b_d)&=&\sum_{\substack{j = 1,\ldots,d:\\ |\ln \omega_j| > \tau b_d}} \omega_j  + R_0(d,\tau b_d),\\
	\sum_{j=1}^{d} \Expec \Bigl[\hat{U}_j\id(|\hat{U}_j|\leqslant \tau b_d)\Bigr]&=& \sum_{\substack{j = 1,\ldots,d:\\ |\ln \omega_j| \leqslant \tau b_d}}\dfrac{ |\ln \omega_j|\,\omega_j}{1 - \omega_j}-R_1(d,\tau b_d), \\
	\sum_{j=1}^{d}\Var \Bigl[\hat{U}_j\id(|\hat{U}_j|\leqslant \tau b_d)\Bigr]&=&  \sum_{\substack{j = 1,\ldots,d:\\ |\ln \omega_j| \leqslant \tau b_d}}\dfrac{ |\ln \omega_j|^2\,\omega_j}{(1 - \omega_j)^2}-R_2(d,\tau b_d),
\end{eqnarray*}
where the functions $R_0$, $R_1$, and $R_2$ are defined as in Lemma \ref{lm_sum_omegaj}. Therefore \eqref{th_nYde_LogAsymp_cond} is equivalent to the following three conditions
\begin{eqnarray*}
		\mathrm{(A')}&&\lim_{d\to\infty}\biggl(\sum_{\substack{j = 1,\ldots,d:\\ |\ln \omega_j| > \tau b_d}} \omega_j  + R_0(d,\tau b_d)\biggr)=-L(\tau)\quad\text{for all}\quad \tau>0;
\end{eqnarray*}
\begin{eqnarray*}
	\mathrm{(B')}&&\lim_{d\to\infty}\dfrac{1}{b_d} \biggl(\sum_{\substack{j = 1,\ldots,d:\\ |\ln \omega_j| \leqslant \tau b_d}}\dfrac{ |\ln \omega_j|\,\omega_j}{1 - \omega_j}-R_1(d,\tau b_d)-\hat{a}_d\biggr)=c+\gamma_\tau \quad\text{for all}\quad \tau>0;
\end{eqnarray*}
\begin{eqnarray*}
	\mathrm{(C')}&&\lim_{\tau\to0}\varliminf_{d\to\infty} \dfrac{1}{b_d^2}\biggl(\sum_{\substack{j = 1,\ldots,d:\\ |\ln \omega_j| \leqslant \tau b_d}}\dfrac{ |\ln \omega_j|^2\,\omega_j}{(1 - \omega_j)^2}-R_2(d,\tau b_d)\biggr)=\\
	&&{}\qquad=\lim_{\tau\to0}\varlimsup_{d\to\infty} \dfrac{1}{b_d^2}\biggl(\sum_{\substack{j = 1,\ldots,d:\\ |\ln \omega_j| \leqslant \tau b_d}}\dfrac{ |\ln \omega_j|^2\,\omega_j}{(1 - \omega_j)^2}-R_2(d,\tau b_d)\biggr)=v^2.
\end{eqnarray*}

We first show that $\mathrm{(A)}$, $\mathrm{(B)}$, $\mathrm{(C)}$ imply $\mathrm{(A')}$, $\mathrm{(B')}$, $\mathrm{(C')}$. Due to the condition $(\mathrm{C})$, there exist $\tau_0>0$ and $C_0>0$  such that for all  $d\in\N$ we have
\begin{eqnarray}\label{ineq_sum_ln2omegaj_tau0}
	\dfrac{1}{b_d^2}\sum_{\substack{j = 1,\ldots,d:\\ |\ln \omega_j| \leqslant \tau_0 b_d}}\dfrac{ |\ln \omega_j|^2\,\omega_j}{(1 - \omega_j)^2}\leqslant C_0.
\end{eqnarray}
Since $\omega_j\in(0,1)$, the inequality $|\ln\omega_j|>1-\omega_j$ always holds. Hence we conclude that
\begin{eqnarray}\label{ineq_sum_omegaj_tau0}
	\dfrac{1}{b_d^2}\sum_{\substack{j = 1,\ldots,d:\\ |\ln \omega_j| \leqslant \tau_0 b_d}}\omega_j\leqslant C_0,\quad d\in\N.
\end{eqnarray}

We first show that   $R_0(d,\tau b_d)\to 0$, $d\to\infty$, for every $\tau>0$. We fix $\tau>0$ and $\tau_*<\min\{\tau_0, \tau\}$. Observe that 
\begin{eqnarray*}
	R_0(d,\tau b_d) = \sum_{\substack{j = 1,\ldots,d:\\ |\ln \omega_j| \leqslant \tau b_d}} \omega_j^{k_{j}(\tau b_d)} \leqslant \sum_{\substack{j = 1,\ldots,d:\\ |\ln \omega_j| \leqslant \tau_* b_d}} \omega_j^{k_{j}(\tau b_d)}+\sum_{\substack{j = 1,\ldots,d:\\ |\ln \omega_j|> \tau_* b_d}} \omega_j^{k_{j}(\tau b_d)}.
\end{eqnarray*}
By the definition of $k_j(\cdot)$ (see Lemma \ref{lm_sum_omegaj}), we have the inequalities $k_j(\tau b_d)\geqslant 2$ and $\omega_j^{k_j(\tau b_d)}\leqslant e^{-\tau b_d}$, which give
\begin{eqnarray*}
	R_0(d,\tau b_d)  \leqslant e^{-\tau b_d}\!\!\! \!\sum_{\substack{j = 1,\ldots,d:\\ |\ln \omega_j| \leqslant \tau_* b_d}} 1+\sum_{\substack{j = 1,\ldots,d:\\ |\ln \omega_j|> \tau_* b_d}} \omega_j^{2}.
\end{eqnarray*}
Using conditions in the sums, we obtain
\begin{eqnarray*}
	R_0(d,\tau b_d)  &\leqslant& e^{-(\tau-\tau_*) b_d}\!\!\!\! \sum_{\substack{j = 1,\ldots,d:\\ |\ln \omega_j| \leqslant \tau_* b_d}} \omega_j+e^{-\tau_* b_d}\!\!\!\!\sum_{\substack{j = 1,\ldots,d:\\ |\ln \omega_j|> \tau_* b_d}} \omega_j\\
	&\leqslant& e^{-(\tau-\tau_*) b_d}\!\!\!\! \sum_{\substack{j = 1,\ldots,d:\\ |\ln \omega_j| \leqslant \tau_0 b_d}} \omega_j+e^{-\tau_* b_d}\!\!\!\!\sum_{\substack{j = 1,\ldots,d:\\ |\ln \omega_j|> \tau_* b_d}} \omega_j.
\end{eqnarray*}
Here the first sum is less than $C_0 b_d^2$ by \eqref{ineq_sum_omegaj_tau0} and the second sum is bounded by some constant $C_1$ due to $\mathrm{(A)}$. Thus 
\begin{eqnarray*}
	R_0(d,\tau b_d)\leqslant C_0 b_d^2 e^{-(\tau-\tau_*) b_d}  + C_1 e^{-\tau_* b_d}.
\end{eqnarray*}
Since $\tau> \tau_*>0$ and $b_d\to \infty$, we obtain that  $R_0(d,\tau b_d)\to 0$, $d\to\infty$. This, together with  $\mathrm{(A)}$, yields $\mathrm{(A')}$.

We now consider 
\begin{eqnarray*}
	R_1(d,\tau b_d)=\sum_{\substack{j = 1,\ldots,d:\\ |\ln \omega_j| \leqslant \tau b_d}} \dfrac{|\ln \omega_j|\,\omega_j^{k_{j}(\tau b_d)}}{1 - \omega_j} \bigl(k_{j}(\tau b_d)(1 - \omega_j) + \omega_j\bigr).
\end{eqnarray*}
According to the definition of $k_j(\cdot)$, we have
\begin{eqnarray}\label{ineq_kjtaubd}
	k_{j}(\tau b_d)\leqslant 2\bigl( k_j(\tau b_d)-1\bigr)\leqslant \dfrac{2\tau b_d }{|\ln \omega_j|}.
\end{eqnarray}
Hence
\begin{eqnarray*}
	R_1(d,\tau b_d)\leqslant	2\tau b_d\sum_{\substack{j = 1,\ldots,d:\\ |\ln \omega_j| \leqslant \tau b_d}} \omega_j^{k_{j}(\tau b_d)}  +  \sum_{\substack{j = 1,\ldots,d:\\ |\ln \omega_j| \leqslant \tau b_d}} \dfrac{|\ln \omega_j|\,\omega_j}{1 - \omega_j}\cdot \omega_j^{k_{j}(\tau b_d)} .
\end{eqnarray*}
Note that the function  $\omega_j\mapsto \tfrac{|\ln \omega_j| \omega_j}{1-\omega_j}$ is bounded by some constant $C_3$ for all values $\omega_j\in(0,1)$. Therefore
\begin{eqnarray}\label{ineq_R1}
	R_1(d,\tau d)\leqslant (2\tau b_d +C_3) \sum_{\substack{j = 1,\ldots,d:\\ |\ln \omega_j| \leqslant \tau b_d}} \omega_j^{k_{j}(\tau b_d)}=(2\tau b_d +C_3)\cdot R_0 (d,\tau b_d).
\end{eqnarray}
Since $R_0(d,\tau b_d)\to 0$, $d\to\infty$, we have $R_1(d,\tau d)=o(b_d)$, $d\to \infty$. Due to $\mathrm{(B)}$, this yields $\mathrm{(B')}$. 

We now consider 
\begin{eqnarray*}
	R_2(d,\tau b_d)&=&\sum_{\substack{j = 1,\ldots,d:\\ |\ln \omega_j| \leqslant \tau b_d}}\dfrac{|\ln \omega_j|^2\omega_j^{k_j(\tau b_d)}}{(1-\omega_j)^2}\biggl(k_j(\tau b_d)^2(1-\omega_j)^2\bigl(1+\omega_j^{k_j(\tau b_d)}\bigr)\\
	&&{}\quad+2k_j(\tau b_d)(1-\omega_j)\omega_j^{k_j(\tau b_d)+1} +\omega_j^{k_j(\tau b_d)+2}+(1-\omega_j)\omega_j\biggr).
\end{eqnarray*}
Using the inequality \eqref{ineq_kjtaubd},
\begin{eqnarray*}
	R_2(d,\tau b_d)&\leqslant&(2\tau b_d)^2\sum_{\substack{j = 1,\ldots,d:\\ |\ln \omega_j| \leqslant \tau b_d}}\omega_j^{k_j(\tau b_d)}\bigl(1+\omega_j^{k_j(\tau b_d)}\bigr)+4\tau b_d\sum_{\substack{j = 1,\ldots,d:\\ |\ln \omega_j| \leqslant \tau b_d}}\dfrac{|\ln \omega_j|}{1-\omega_j}\omega_j^{2k_j(\tau b_d)+1}\\
	&&{}+\sum_{\substack{j = 1,\ldots,d:\\ |\ln \omega_j| \leqslant \tau b_d}}\dfrac{|\ln \omega_j|^2\omega_j^{2k_j(\tau b_d)+2}}{(1-\omega_j)^2}+\sum_{\substack{j = 1,\ldots,d:\\ |\ln \omega_j| \leqslant \tau b_d}}\dfrac{|\ln \omega_j|^2\omega_j^{k_j(\tau b_d)+1}}{1-\omega_j}.
\end{eqnarray*}
The function  $\omega_j\mapsto \tfrac{|\ln \omega_j|^2 \omega_j}{1-\omega_j}$ is bounded by some constant $C_4$ for all values $\omega_j\in(0,1)$. Using this and the same fact about $\omega_j\mapsto \tfrac{|\ln \omega_j| \omega_j}{1-\omega_j}$, we obtain
\begin{eqnarray*}
	R_2(d,\tau b_d)&\leqslant&4\tau^2 b_d^2\sum_{\substack{j = 1,\ldots,d:\\ |\ln \omega_j| \leqslant \tau b_d}}\omega_j^{k_j(\tau b_d)}\bigl(1+\omega_j^{k_j(\tau b_d)}\bigr)+4\tau b_d\sum_{\substack{j = 1,\ldots,d:\\ |\ln \omega_j| \leqslant \tau b_d}}\omega_j^{2k_j(\tau b_d)}\\
	&&{}+C_3^2\sum_{\substack{j = 1,\ldots,d:\\ |\ln \omega_j| \leqslant \tau b_d}}\omega_j^{k_j(\tau b_d)}+C_4\sum_{\substack{j = 1,\ldots,d:\\ |\ln \omega_j| \leqslant \tau b_d}}\omega_j^{2k_j(\tau b_d)}.
\end{eqnarray*}
Next, since $\omega_j^{k_j(\tau b_d)}\leqslant 1$, we have
\begin{eqnarray}
	R_2(d,\tau b_d)&\leqslant&(8\tau^2 b_d^2+4\tau b_d+C_3^2 +C_4)\sum_{\substack{j = 1,\ldots,d:\\ |\ln \omega_j| \leqslant \tau b_d}}\omega_j^{k_j(\tau b_d)}\nonumber\\
	&=&(8\tau^2 b_d^2+4\tau b_d+C_3^2 +C_4) R_0(d,\tau b_d).\label{ineq_R2}
\end{eqnarray}
Here $R_0(d,\tau b_d)\to 0$, $d\to\infty$. Hence $R_2(d,\tau b_d)=o(b_d^2)$, $d\to \infty$. Due to $\mathrm{(C)}$, this yields $\mathrm{(C')}$. 

We now show that $\mathrm{(A')}$, $\mathrm{(B')}$, $\mathrm{(C')}$ imply $\mathrm{(A)}$, $\mathrm{(B)}$, $\mathrm{(C)}$. From $\mathrm{(A')}$ it follows  that  $\sup_{d \in\N} R_0(d,\tau b_d)<\infty$  for every $\tau>0$. By this relation and \eqref{ineq_R2}, we have
\begin{eqnarray*}
	\sup_{d \in\N}\dfrac{1}{b_d^2}\,  R_2(d,\tau b_d)<\infty\quad \text{for every  }\tau>0.
\end{eqnarray*}
This and $\mathrm{(C')}$ yield \eqref{ineq_sum_ln2omegaj_tau0} for every $d\in\N$ and some $\tau_0>0$ and $C_0>0$. Hence, as we showed above, it follows that  $R_0(d,\tau b_d)\to 0$, $d\to\infty$, for every $\tau>0$. Here we use
\begin{eqnarray*}
	\sup_{d \in\N} \sum_{\substack{j = 1,\ldots,d:\\ |\ln \omega_j| > \tau b_d}}\omega_j<\infty\quad \text{for every } \tau>0,
\end{eqnarray*}
which  follows from $\mathrm{(A')}$. Next, from \eqref{ineq_R1} and \eqref{ineq_R2} we correspondingly obtain that $R_1(d,\tau b_d)=o(b_d)$ and $R_2(d,\tau b_d)=o(b_d^2)$, $d\to\infty$, for every $\tau>0$. These relations for $R_k(d,\tau b_d)$, $k=1,2,3$, and conditions $\mathrm{(A')}$, $\mathrm{(B')}$, $\mathrm{(C')}$ give $\mathrm{(A)}$, $\mathrm{(B)}$, $\mathrm{(C)}$. \quad $\Box$

\section{Examples}
In this section we will study asymptotics for $n^{Y_d}(\e)$ under particular additional assumptions on the sequence of the length scale parameters $\sigma_j$, $j\in\N$. 

We begin with the case, when the
sequence $(\sigma_j)_{j \in \N}$ tends to to a non-negative constant $\sigma$. Observe that here approximation complexity $n^{Y_d}(\e)$ is unbounded due to Proposition \ref{pr_unboundedness}. The following proposition specifies asymptotic \eqref{th_nYde_LogAsymp_cond} for this case.

\begin{Proposition}	\label{pr_non_zero}
	Suppose that $\sigma_j \to \sigma$, $0 \leqslant \sigma < \infty$.
	Then
      \begin{eqnarray}	\label{eq_non_zero}
	\forall \e\in(0,1)\quad	\ln n^{Y_d}(\e) = a_d + \Phi^{-1}(1 - \e^2) b_d + o(b_d), \quad d \to \infty,  
	\end{eqnarray}
	where $\Phi$ is the  distribution function of the standard Gaussian law,
	\begin{eqnarray*}
		a_d = \sum\limits_{j = 1}^d \biggl( \frac{|\ln \omega_j| \omega_j}{1 - \omega_j} + |\ln (1 - \omega_j)| \biggr),\qquad b_d = \rho \sqrt{d},\quad d \in \N,
	\end{eqnarray*}
\begin{eqnarray*}
 \rho=\lim\limits_{j \to \infty} \tfrac{|\ln \omega_j| \sqrt{\omega_j}}{1 - \omega_j}=
 \begin{cases}
 	\tfrac{|\ln \omega| \sqrt{\omega}}{1 - \omega},& \sigma>0,\\
 	1,& \sigma=0,\\
 \end{cases}\qquad \omega=\lim\limits_{j \to \infty} \omega_j= \biggr(1+\tfrac{\sigma^2}{2} \Bigl(1+\sqrt{1+\tfrac{4}{\sigma^2}}\,\Bigr)\biggr)^{-1}\text{ for}\quad\sigma>0.
\end{eqnarray*}
\end{Proposition}

\begin{Remark}
	If $\sigma>0$ then 
	\begin{eqnarray*}
		a_d= \biggl( \dfrac{|\ln \omega| \omega}{1 - \omega} + |\ln (1 - \omega)| \biggr) d +o(d),\quad d\to\infty.
	\end{eqnarray*}
\end{Remark}
Indeed, since $\lim\limits_{j \to \infty} \omega_j=\omega$ we have
\begin{eqnarray*}
	a_d\sim d\cdot \lim\limits_{d\to\infty}\biggl( \frac{|\ln \omega_d| \omega_d}{1 - \omega_d} + |\ln (1 - \omega_d)| \biggr)=d\cdot \biggl( \dfrac{|\ln \omega| \omega}{1 - \omega} + |\ln (1 - \omega)| \biggr),\quad d\to\infty.\\
\end{eqnarray*}

\begin{Remark}
	If $\sigma=0$ then $a_d$ can increase arbitrarily fast.
\end{Remark}
It is seen from the inequality 
\begin{eqnarray*}
	a_d\geqslant |\ln(1-\omega_d)|,\quad d\in\N,
\end{eqnarray*}
and from the relation $\lim_{d \to \infty} \omega_d=1$,  which holds under the assumption $\sigma=0$.\\

More sharp asymptotics of $a_d$ (up to $o(b_d)$) require additional assumptions on the rate of convergence of $\sigma_j$ to $\sigma$ as $j\to\infty$.\\

\noindent
\textbf{Proof of Proposition \ref{pr_non_zero}.}\quad It is well known that $\Phi$ is a self-decomposable distribution function with the triplet $(0, 1, L)$, where $L(x) = 0$ for all $x \in \R \setminus \{0\}$. By Theorem \ref{th_nYde_LogAsymp} it is sufficient to check the following ensemble of conditions to obtain the asymptotics \eqref{eq_non_zero}:
\begin{eqnarray}
	&&\lim_{d \to \infty} \sum_{\substack{j = 1,\ldots, d\\ |\ln \omega_j| > \tau b_d}} \omega_j = 0,\quad \tau>0;\label{cond_tails_normal}\\
	&&\lim_{d \to \infty}  \frac{1}{b_d} \biggl(\sum_{\substack{j = 1,\ldots, d\\ |\ln \omega_j| \leqslant \tau b_d}} \frac{ |\ln \omega_j|\omega_j}{1 - \omega_j} - \hat{a}_d \biggr) =
	0,\quad \tau > 0; \label{cond_mean_normal}\\
	&&\lim_{\tau \to 0} \varliminf_{d \to \infty} \frac{1}{b_d^2} \sum_{\substack{j = 1,\ldots, d\\ |\ln \omega_j| \leqslant \tau b_d}} \frac{ |\ln \omega_j|^2\omega_j}{(1 - \omega_j)^2} =
	\lim_{\tau \to 0} \varlimsup_{d \to \infty} \frac{1}{b_d^2} \sum_{\substack{j = 1,\ldots, d\\ |\ln \omega_j| \leqslant \tau b_d}} \frac{ |\ln \omega_j|^2\omega_j}{(1 - \omega_j)^2}  = 1.\label{cond_var_normal} 
\end{eqnarray}
Here $\hat{a}_d$ is defined by \eqref{def_hatad}.

Observe that $\lim_{j\to\infty}\omega_j \in (0, 1]$, and, consequently, the sequence $(|\ln \omega_j|)_{j \in \N}$ is bounded. Since $b_d \to \infty$, $d \to \infty$, for every $\tau>0$ there exists $d_\tau\in\N$ such that $|\ln \omega_j| \leqslant \tau b_d$ holds for any $j\in\N$ and $d\geqslant d_\tau$.
Therefore for every $\tau>0$ we have
\begin{eqnarray*}
	\sum_{\substack{j = 1,\ldots, d\\ |\ln \omega_j| > \tau b_d}} \omega_j = 0,\quad\text{and}\quad 	\sum_{\substack{j = 1,\ldots, d\\ |\ln \omega_j| \leqslant \tau b_d}} \dfrac{ |\ln \omega_j|\omega_j}{1 - \omega_j} - \hat{a}_d =\sum_{j = 1}^d \dfrac{ |\ln \omega_j|\omega_j}{1 - \omega_j} - \hat{a}_d = 0
\end{eqnarray*}
for all sufficiently large $d$, i.e.  we get $\eqref{cond_tails_normal}$ and $\eqref{cond_mean_normal}$. Next, we also have
\begin{eqnarray*}
	\frac{1}{b_d^2} \sum_{\substack{j = 1,\ldots, d\\ |\ln \omega_j| \leqslant \tau b_d}} \frac{ |\ln \omega_j|^2\omega_j}{(1 - \omega_j)^2} = \frac{1}{b_d^2} \sum_{j=1}^d \dfrac{ |\ln \omega_j|^2\omega_j}{(1 - \omega_j)^2}=\dfrac{1}{\rho^2}\cdot
	\frac{1}{d}\sum_{j = 1}^d  \dfrac{ |\ln \omega_j|^2\omega_j}{(1 - \omega_j)^2}
\end{eqnarray*}
for all sufficiently large $d$.  The arithmetic mean in the right-hand side tends to $\lim_{j\to\infty}\tfrac{ |\ln \omega_j|^2\omega_j}{(1 - \omega_j)^2}=\rho^2$. Thus for every $\tau>0$ we get
\begin{eqnarray*}
	\lim_{d\to\infty} \frac{1}{b_d^2} \sum_{\substack{j = 1,\ldots, d\\ |\ln \omega_j| \leqslant \tau b_d}} \frac{ |\ln \omega_j|^2\omega_j}{(1 - \omega_j)^2}=1,
\end{eqnarray*}
that implies $\eqref{cond_var_normal}$.\quad $\Box$\\

We now proceed to the case  $\sigma_j \to \infty$ as $j \to \infty$. Here the asymptotics of $n^{Y_d}(\e)$ depends on the velocity of $\sigma_j$. In order to evidently illustrate the application of the general results from the previous section and to avoid routine unwieldy calculations, we will assume that
\begin{eqnarray}\label{assum_sigmaj}
	\sigma_j^2 \sim \beta j^{\alpha},\quad j\to\infty, 
\end{eqnarray}
with some $\alpha>0$ and $\beta>0$.

The following assertion directly follows from Propositions \ref{pr_boundedness} and \ref{pr_unboundedness}.

\begin{Proposition}
	Let $(\sigma_j)_{j\in\N}$ satisfy \eqref{assum_sigmaj} with some $\alpha>0$ and $\beta>0$. If $\alpha>1$ then $\sup_{d \in\N} n^{Y_d}(\e)<\infty$ for every $\e\in(0,1)$. If $\alpha\leqslant1$ then $n^{Y_d}(\e)\to\infty$ as $d\to\infty$ for every $\e\in(0,1)$. 
\end{Proposition}

We now propose the asymptotics of $n^{Y_d}(\e)$ for the case $\alpha\in(0,1]$.

\begin{Proposition}	\label{pr_zero_normal_dist}
	Let $(\sigma_j)_{j\in\N}$ satisfy \eqref{assum_sigmaj} with some $\alpha\in(0,1]$ and $\beta>0$. Then
	\begin{eqnarray*}
		\forall \e\in(0,1)\quad\ln n^{Y_d}(\e) = a_d + \Phi^{-1}(1 - \e^2) b_d + o(b_d), \quad d \to \infty,  
	\end{eqnarray*}
	where $\Phi$ is the  distribution function of the standard Gaussian law,
	\begin{eqnarray*}
		a_d = \sum\limits_{j = 1}^d \biggl( \frac{|\ln \omega_j| \omega_j}{1 - \omega_j} + |\ln (1 - \omega_j)| \biggr),\quad\text{and}\quad 
			b_d = 
			\begin{cases}
				\tfrac{ \alpha  }{\sqrt{(1 - \alpha)\beta}}\,d^{\tfrac{1 -  \alpha}{2}} \ln d,&  \alpha \in (0, 1),\\
				\tfrac{1}{\sqrt{3\beta}}\,(\ln d)^{3 / 2},& \alpha  = 1,
			\end{cases}
		\quad d\in\N.
		\end{eqnarray*}
\end{Proposition}
\begin{Remark}
	The following asymptotic holds
	\begin{eqnarray}\label{rem_ad}
		a_d\sim 
		\begin{cases}
			\tfrac{\alpha}{(1-\alpha)\beta}\,d^{1-\alpha}\ln d+\tfrac{1}{(1-\alpha)\beta}\,d^{1-\alpha},& \alpha \in (0, 1),\\
			\tfrac{1}{2\beta}\,(\ln d)^2+ \tfrac{1}{\beta}\,\ln d,& \alpha=1,
		\end{cases}
		\quad d\to\infty.
	\end{eqnarray}
\end{Remark}
Indeed, since $\sigma_j \to \infty$ as $j \to \infty$, we have $\omega_j \sim \sigma_j^{-2}\sim (\beta j^{\alpha})^{-1}$, $j\to\infty$,  due to \eqref{def_omegaj} and \eqref{assum_sigmaj}.  Then
\begin{eqnarray*}
	a_d\sim \sum_{j=1}^d\biggl( \dfrac{\alpha\ln j}{ \beta j^{\alpha}} +  \dfrac{1}{\beta j^{\alpha}}\biggr)\sim \dfrac{\alpha}{\beta} \int\limits_{1}^{d} \dfrac{\ln x}{x^\alpha}\, \dd x + \dfrac{1}{\beta} \int\limits_{1}^{d} \dfrac{1}{x^\alpha}\, \dd x,\quad d\to\infty.
\end{eqnarray*}
Using known asymptotic relations (see \cite{BingGoldTeug} and \cite{MakGolLodPod}), we obtain \eqref{rem_ad}.\\

More sharp asymptotics of $a_d$ (up to $o(b_d)$) require additional assumptions on the asymptotics of the difference of $\sigma_j^2- \beta j^{\alpha}$ as $j\to\infty$.\\

\noindent
\textbf{Proof of Proposition \ref{pr_zero_normal_dist}.}\quad As in Proposition \ref{pr_non_zero}, according to Theorem \ref{th_nYde_LogAsymp}, it is sufficient to check  conditions \eqref{cond_tails_normal}--\eqref{cond_var_normal} with given $a_d$, $b_d$, and with $\hat{a}_d$ defined by \eqref{def_hatad}, $d\in\N$. 

Since $\omega_j \sim (\beta j^{\alpha})^{-1}$, $j\to\infty$, it is not difficult to see that for every $\tau>0$
\begin{eqnarray*}
	\dfrac{1}{\tau b_d}\max_{j\in\{1,\ldots, d\}} |\ln \omega_j|\to 0,\quad d\to\infty.
\end{eqnarray*}
Therefore for every $\tau>0$ we have
\begin{eqnarray*}
	\sum_{\substack{j = 1,\ldots, d\\ |\ln \omega_j| > \tau b_d}} \omega_j = 0,\quad\text{and}\quad 	\sum_{\substack{j = 1,\ldots, d\\ |\ln \omega_j| \leqslant \tau b_d}} \dfrac{ |\ln \omega_j|\omega_j}{1 - \omega_j} - \hat{a}_d =\sum_{j = 1}^d \dfrac{ |\ln \omega_j|\omega_j}{1 - \omega_j} - \hat{a}_d = 0
\end{eqnarray*}
for all sufficiently large $d$, i.e.  we get \eqref{cond_tails_normal} and \eqref{cond_mean_normal}. Next, we also have
\begin{eqnarray*}
	\sum_{\substack{j = 1,\ldots, d\\ |\ln \omega_j| \leqslant \tau b_d}} \frac{ |\ln \omega_j|^2\omega_j}{(1 - \omega_j)^2} = \sum_{j=1}^d \dfrac{ |\ln \omega_j|^2\omega_j}{(1 - \omega_j)^2}
\end{eqnarray*}
for every $\tau>0$ and for all sufficiently large $d$.  Using the asimptotics $\omega_j \sim (\beta j^{\alpha})^{-1}$, $j\to\infty$, we get
\begin{eqnarray*}
	\sum_{j=1}^d \dfrac{ |\ln \omega_j|^2\omega_j}{(1 - \omega_j)^2}=	\sum_{j=1}^d \dfrac{\bigl(\alpha\ln j+\ln \beta+o(1)\bigr)^2}{ \beta j^{\alpha}(1+o(1))}\sim \dfrac{\alpha^2}{\beta}	\sum_{j=1}^d \dfrac{  (\ln j)^2}{j^{\alpha}},\quad d\to\infty.
\end{eqnarray*}
By known asymptotic relations (see \cite{BingGoldTeug} and \cite{MakGolLodPod}),  we obtain
\begin{eqnarray*}
	\sum_{j=1}^d \dfrac{  (\ln j)^2}{j^{\alpha}}\sim \int\limits_1^d \dfrac{ (\ln x)^2}{x^{\alpha}}\, \dd x\sim
	\begin{cases}
		\tfrac{1}{1-\alpha}\,d^{1-\alpha}(\ln d)^2,& \alpha \in (0, 1),\\
		\tfrac{1}{3}(\ln d)^3,& \alpha=1,
	\end{cases}
\quad d\to\infty,
\end{eqnarray*}
Hence for every $\tau>0$
\begin{eqnarray*}
	\sum_{\substack{j = 1,\ldots, d\\ |\ln \omega_j| \leqslant \tau b_d}} \frac{ |\ln \omega_j|^2\omega_j}{(1 - \omega_j)^2}\sim
	\begin{cases}
		\tfrac{\alpha^2}{(1-\alpha)\beta}\,d^{1-\alpha}(\ln d)^2,& \alpha \in (0, 1),\\
		\tfrac{1}{3\beta}(\ln d)^3,& \alpha=1,
	\end{cases}
	\quad d\to\infty.
\end{eqnarray*}
Thus
\begin{eqnarray*}
	\lim_{d\to\infty} \frac{1}{b_d^2} \sum_{\substack{j = 1,\ldots, d\\ |\ln \omega_j| \leqslant \tau b_d}} \frac{ |\ln \omega_j|^2\omega_j}{(1 - \omega_j)^2}=1,\quad \tau>0,
\end{eqnarray*}
that implies $\eqref{cond_var_normal}$.\quad $\Box$\\

There exist particular cases, where the $\e$-component of asymptotics \eqref{th_nYde_LogAsymp_cond} is a quantile of a non-Gaussian self-decomposable distribution function. One such case is considered in the next proposition. There appear  distribution functions $D_\mu$, $\mu>0$,  of $\mu$-convolution powers of the Dickman law (see \cite{Hensl} and \cite{Khart1}). It is known that $D_\mu$ has the spectral representation triplet $(\tfrac{\pi\mu }{4}, 0, L)$, where $L(x)=\mu \ln x\, \id(x\in(0,1])$, $x\in\R\setminus\{0\}$.

\begin{Proposition}	\label{pr_dickman}
	Suppose that $\sigma_j^2 \sim  \beta j \ln j$,  $j \to \infty$, with some $\beta > 0$. Then 
      \begin{eqnarray*}
		\forall \e\in(0,1)\quad \ln n^{Y_d}(\e) = D^{-1}_{1/\beta}(1 - \e^2) \ln d + o(\ln d), \quad d \to \infty.  
	\end{eqnarray*}
\end{Proposition}

\noindent
\textbf{Proof of Proposition \ref{pr_dickman}.}\quad We use Theorem \ref{th_nYde_LogAsymp} with  $a_d\defeq 0$, $b_d\defeq \max\{\ln d,1\}$, $d\in\N$, $G\defeq D_{1/\beta}$, $c= \tfrac{\pi}{4\beta}$, $v=0$, and $L(x)=\tfrac{1}{\beta}\ln x\, \id(x\in(0,1])$, $x\in\R\setminus\{0\}$.  So it is sufficient to check the following ensemble of conditions:
\begin{eqnarray}
	&&\lim_{d \to \infty} \sum_{\substack{j = 1,\ldots, d\\ |\ln \omega_j| > \tau \ln d}} \omega_j = -\tfrac{1}{\beta} \ln \tau \id(\tau \in (0, 1]),\quad \tau>0,\label{pr_dickman_cond_tails}\\
	&&\lim_{d \to \infty}  \frac{1}{\ln d} \biggl(\sum_{\substack{j = 1,\ldots, d\\ |\ln \omega_j| \leqslant \tau \ln d}} \frac{ |\ln \omega_j|\omega_j}{1 - \omega_j}+\sum\limits_{j=1}^{d} |\ln(1-\omega_j)|\biggr) =
	\frac{\pi}{4\beta} + \gamma_\tau,\quad \tau > 0, \label{pr_dickman_cond_mean}\\
	&&\lim_{\tau \to 0} \varliminf_{d \to \infty} \frac{1}{(\ln d)^2}\!\!\! \sum_{\substack{j = 1,\ldots, d\\ |\ln \omega_j| \leqslant \tau \ln d}} \frac{ |\ln \omega_j|^2\omega_j}{(1 - \omega_j)^2} =
	\lim_{\tau \to 0} \varlimsup_{d \to \infty} \frac{1}{(\ln d)^2}\!\!\! \sum_{\substack{j = 1,\ldots, d\\ |\ln \omega_j| \leqslant \tau \ln d}} \frac{ |\ln \omega_j|^2\omega_j}{(1 - \omega_j)^2}  = 0.\label{pr_dickman_cond_var}
\end{eqnarray}
Here
\begin{eqnarray*}
	\gamma_\tau = \dfrac{1}{\beta}\int\limits_0^\tau \dfrac{x^2 \id(x \in (0, 1))}{1 + x^2}\, \dd x - \dfrac{1}{\beta} \int\limits_\tau^{+\infty}\frac{\id(x \in (0, 1))}{1 + x^2}\, \dd x,\quad \tau>0.
\end{eqnarray*}
Due to the assumption, we have 
\begin{eqnarray} \label{pr_dickman_cond_omegaj}
	\omega_j \sim \sigma_j^{-2}\sim \dfrac{1}{\beta j \ln j}, \quad j \to \infty.
\end{eqnarray}
We first check \eqref{pr_dickman_cond_tails}. It is easily seen that for $\tau>1$ we have $|\ln \omega_j|\leqslant \tau \ln d$, $j=1,\ldots, d$, for all sufficiently large $d$. Therefore \eqref{pr_dickman_cond_tails} obviously holds in this case. For $\tau\in (0,1]$ we set
\begin{eqnarray}\label{def_jdtau}
	j_{d,\tau}=\min\{j\in\N: |\ln \omega_j|>\tau \ln d\},\quad  d\in\N.
\end{eqnarray}
Due to \eqref{pr_dickman_cond_omegaj}, it is not difficult to see that $j_{d,\tau}\leqslant d$ for all sufficiently large $d$ and $\ln j_{d,\tau}\sim \tau |\ln d|$, $d\to\infty$. Next, observe that
\begin{eqnarray*}
	\sum_{\substack{j = 1,\ldots, d\\ |\ln \omega_j| > \tau b_d}}\omega_j=
	\sum_{j = j_{d,\tau}}^d\omega_j\sim \sum_{j = j_{d,\tau}}^d \dfrac{1}{\beta j \ln j},\quad d\to\infty.
\end{eqnarray*}
Using the known asymptotics (see \cite{MakGolLodPod}, 2.13, p. 21)
\begin{eqnarray}\label{asymp_klnk}
	\sum_{k=1}^{n} \dfrac{1}{k \ln k}= \ln \ln n + C + o(1),\quad n\to\infty,
\end{eqnarray}
with a constant $C$, we obtain
\begin{eqnarray*}
	\sum_{j = j_{d,\tau}}^d \dfrac{1}{\beta j \ln j}&=& \tfrac{1}{\beta} (\ln\ln d - \ln\ln j_{d,\tau})+o(1)\\
	&=&\tfrac{1}{\beta} (\ln\ln d - \ln\ln d^\tau)+o(1)\\
	&=&-\tfrac{1}{\beta} \ln \tau + o(1),\quad d\to\infty.
\end{eqnarray*}
Thus \eqref{pr_dickman_cond_tails} holds.

We now check \eqref{pr_dickman_cond_mean}. Observe that for every $\tau>0$
\begin{eqnarray}
    \frac{\pi}{4\beta} + \gamma_\tau &=& \frac{\pi}{4\beta} +
    \dfrac{1}{\beta}\int\limits_0^{\min\{\tau, 1\}} \frac{x^2}{1 + x^2}\, \dd x - \dfrac{1}{\beta} \int\limits_{\min\{\tau, 1\}}^{1}\frac{1}{1 + x^2}\, \dd x \nonumber\\
     &=& \frac{\pi}{4\beta} +
    \dfrac{1}{\beta}\int\limits_0^{\min\{\tau, 1\}} \dd x - \dfrac{1}{\beta} \int\limits_{0}^{1}\frac{1}{1 + x^2}\, \dd x\nonumber\\
    &=&\tfrac{1}{\beta} \min\{\tau,1\}. \label{conc_mintau1}
\end{eqnarray}
Next, due to \eqref{pr_dickman_cond_omegaj} and \eqref{asymp_klnk}, we have
\begin{eqnarray}\label{conc_sum_ln1minusomegaj}
	\sum_{j = 1}^d |\ln(1 - \omega_j)| \sim \sum_{j = 1}^d \omega_j \sim \sum_{j = 1}^d \dfrac{1}{\beta j \ln j} \sim \tfrac{1}{\beta} \ln \ln d = o(\ln d),\quad d \to \infty.
\end{eqnarray}
Next, if $\tau>1$ then $|\ln \omega_j|\leqslant \tau \ln d$, $j=1,\ldots, d$, for all sufficiently large $d$. Therefore
\begin{eqnarray*}
	\sum_{\substack{j = 1,\ldots, d\\ |\ln \omega_j| \leqslant \tau b_d}} \frac{ |\ln \omega_j|\omega_j}{1 - \omega_j}= \sum_{j=1}^d \frac{ |\ln \omega_j|\omega_j}{1 - \omega_j}\sim \sum\limits_{j = 1}^{d}\frac{1}{\beta j} \sim \tfrac{1}{\beta}\, \ln d,\quad  d \to \infty.
\end{eqnarray*}
If $\tau\in(0,1]$ then, according to \eqref{def_jdtau}, we write
\begin{eqnarray*}
	\sum_{\substack{j = 1,\ldots, d\\ |\ln \omega_j| \leqslant \tau b_d}} \frac{ |\ln \omega_j|\omega_j}{1 - \omega_j}= \sum_{j=1}^{j_{d,\tau}-1} \frac{ |\ln \omega_j|\omega_j}{1 - \omega_j}\sim \sum_{j=1}^{j_{d,\tau}-1}\frac{1}{\beta j} \sim \tfrac{1}{\beta}\, \ln j_{d,\tau}\sim \tfrac{\tau}{\beta}\, \ln d,\quad  d \to \infty.
\end{eqnarray*}
Therefore we have for every $\tau>0$
\begin{eqnarray}\label{conc_dickman_means}
	\sum_{\substack{j = 1,\ldots, d\\ |\ln \omega_j| \leqslant \tau b_d}} \frac{ |\ln \omega_j|\omega_j}{1 - \omega_j}\sim \tfrac{1}{\beta}\,\min\{\tau,1\} \ln d,\quad  d \to \infty.
\end{eqnarray}
Thus \eqref{conc_mintau1}, \eqref{conc_sum_ln1minusomegaj}, and \eqref{conc_dickman_means} together yield \eqref{pr_dickman_cond_mean}.

Next, we check the condition \eqref{pr_dickman_cond_var}. Let us fix $\tau\in(0,1)$ and consider
\begin{eqnarray*}
	 \sum_{\substack{j = 1,\ldots, d\\ |\ln \omega_j| \leqslant \tau b_d}} \frac{ |\ln \omega_j|^2\omega_j}{(1 - \omega_j)^2} 
	=  \sum_{j = 1}^{j_{d,\tau}-1} \frac{ |\ln \omega_j|^2\omega_j}{(1 - \omega_j)^2}\sim \sum_{j = 1}^{j_{d,\tau}-1}  \dfrac{\ln j}{\beta j}\sim  \dfrac{(\ln j_{d,\tau})^2}{2\beta}\sim  \dfrac{\tau^2(\ln d)^2}{2\beta},\quad d\to \infty.
\end{eqnarray*}
Therefore
\begin{eqnarray*}
	 \lim_{d \to \infty} \frac{1}{(\ln d)^2} \sum_{\substack{j = 1,\ldots, d\\ |\ln \omega_j| \leqslant \tau b_d}} \frac{ |\ln \omega_j|^2\omega_j}{(1 - \omega_j)^2} =
\frac{\tau^2}{2\beta}.
\end{eqnarray*}
From this we obtain \eqref{pr_dickman_cond_var}.\quad $\Box$\\

\section{Acknowledgments}
The work of A. A. Khartov was supported by RFBR--DFG grant 20-51-12004.

\end{document}